\documentclass[a4paper,12pt]{article}
\usepackage{latexsym}
\usepackage{amsmath}
\usepackage{amssymb}
\usepackage{amscd}
\usepackage{array}
\usepackage[all]{xy}
\usepackage{amsthm}
\usepackage{amsfonts}
\usepackage{enumerate}
\usepackage{multirow}

\usepackage{pdflscape}
\newtheorem{theorem}{Theorem}[]
\newtheorem{proposition}[theorem]{Proposition}

\newtheorem{corollary}[theorem]{Corollary}

\theoremstyle{definition}

\newcommand{\HH}{\mathcal H}

\newcommand{\Z}{\mathbb Z}

\newcommand{\F}{\mathbb F}

\newcommand{\Gal}{\mathrm{Gal}}

\newcommand{\Hol}{\mathrm{Hol}}
\newcommand{\Norm}{\mathrm{Norm}}
\newcommand{\Sym}{\operatorname{Sym}}

\newcommand{\Stab}{\operatorname{Stab}}

\newcommand{\Perm}{\mathrm{Perm}}

\newcommand{\GL}{\mathrm{GL}}

\newcommand{\End}{\operatorname{End}}

\newcommand{\Aut}{\operatorname{Aut}}

\newcommand{\Id}{\operatorname{Id}}
\newcommand{\Syl}{\operatorname{Syl}}

\newcommand{\wL} {{\widetilde{L}}}

\begin{document}
\begin{center}

\Large
Hopf Galois structures on separable field extensions of odd prime power degree

\vspace{1cm}
\large
Teresa Crespo and Marta Salguero

\vspace{0.3cm}
\footnotesize

Departament de Matem\`atiques i Inform\`atica, Universitat de Barcelona (UB), Gran Via de les
Corts Catalanes 585, E-08007 Barcelona, Spain, e-mail: teresa.crespo@ub.edu, msalguga11@alumnes.ub.edu

\end{center}

\date{\today}

\let\thefootnote\relax\footnotetext{Both authors acknowledge support by grant MTM2015-66716-P (MINECO/FEDER, UE).}

\begin{abstract} A Hopf Galois structure on a finite field extension $L/K$ is a pair $(\HH,\mu)$, where $\HH$ is a finite cocommutative $K$-Hopf algebra and $\mu$ a Hopf action. In this paper, we present several results on Hopf Galois structures on odd prime power degree separable field extensions. We prove that if a separable field extension of odd prime power degree has a Hopf Galois structure of cyclic type, then it has no structure of noncyclic type. We determine the number of Hopf Galois structures of cyclic type on a separable field extension of degree $p^n$, $p$ an odd prime, such that the Galois group of its normal closure is a semidirect product $C_{p^n}\rtimes C_D$ of the cyclic group of order $p^n$ and a cyclic group of order $D$, with $D$ prime to $p$.
We characterize the transitive groups of degree $p^3$ which are Galois groups of the normal closure of a separable field extension having some cyclic Hopf Galois structure and determine the number of those. We prove that if a separable field extension of degree $p^3$ has a nonabelian Hopf Galois structure then it has an abelian structure whose type has the same exponent as the nonabelian type. We obtain that, for $p>3$, the two abelian noncyclic Hopf Galois structures do not occur on the same separable extension of degree $p^3$. We present a table which gives the number of Hopf Galois structures of each possible type on a separable extension of degree $27$ to illustrate that for $p=3$, all four noncyclic Hopf Galois structures may occur on the same extension. Finally, putting together all previous results, we list all possible sets of Hopf Galois structure types on a separable extension of degree $p^3$, for $p>3$ a prime.

\end{abstract}

\normalsize

\section{Introduction}
A Hopf Galois structure on a finite extension of fields $L/K$ is a pair $(\HH,\mu)$, where $\HH$ is
a finite cocommutative $K$-Hopf algebra  and $\mu$ is a
Hopf action of $\HH$ on $L$, i.e a $K$-linear map $\mu: \HH \to
\End_K(L)$ giving $L$ a left $\HH$-module algebra structure and inducing a $K$-vector space isomorphism $L\otimes_K \HH\to\End_K(L)$.
Hopf Galois structures were introduced by Chase and Sweedler in \cite{C-S}.
For separable field extensions, Greither and
Pareigis \cite{G-P} give the following group-theoretic
equivalent condition to the existence of a Hopf Galois structure.

\begin{theorem}\label{G-P}
Let $L/K$ be a separable field extension of degree $g$, $\wL$ its Galois closure, $G=\Gal(\wL/K), G'=\Gal(\wL/L)$. Then there is a bijective correspondence
between the set of isomorphism classes of Hopf Galois structures on $L/K$ and the set of
regular subgroups $N$ of the symmetric group $S_g$ normalized by $\lambda (G)$, where
$\lambda:G \hookrightarrow S_g$ is the monomorphism given by the action of
$G$ on the left cosets $G/G'$.
\end{theorem}

For a given Hopf Galois structure on a separable field extension $L/K$ of degree $g$, we will refer to the isomorphism class of the corresponding group $N$ as the type of the Hopf Galois
structure. The Hopf algebra $\HH$ corresponding to a regular subgroup $N$ of $S_g$ normalized by $\lambda (G)$ is the sub-$K$-Hopf algebra $\wL[N]^G$ of the group algebra $\wL[N]$ fixed under the action of $G$, where $G$ acts on $\wL$ by $K$-automorphisms and on $N$ by conjugation through $\lambda$. The Hopf action is induced by $n \mapsto n^{-1}(\overline{1})$, for $n \in N$, where we identify $S_g$ with the group of permutations of $G/G'$ and $\overline{1}$ denotes the class of $1_G$ in $G/G'$.

Childs \cite{Ch1} gives an equivalent  condition to the existence of a Hopf Galois structure introducing the holomorph of the regular subgroup $N$ of $S_g$. We state the more precise formulation of this result due to Byott \cite{B} (see also \cite{Ch2} Theorem 7.3).

\begin{theorem}\label{theoB} Let $G$ be a finite group, $G'\subset G$ a subgroup and $\lambda:G\to \Sym(G/G')$ the morphism given by the action of
$G$ on the left cosets $G/G'$.
Let $N$ be a group of
order $[G:G']$ with identity element $e_N$. Then there is a
bijection between the two sets
$$
{\cal N}=\{\alpha:N\hookrightarrow \Sym(G/G') \mbox{ such that
}\alpha (N)\mbox{ is regular}\}
$$
and
$$
{\cal G}=\{\beta:G\hookrightarrow \Sym(N) \mbox{ such that }\beta
(G')\mbox{ is the stabilizer of } e_N\}
$$

Under this bijection, if $\alpha, \alpha' \in {\cal N}$ correspond to
$\beta, \beta' \in {\cal G}$, then
\begin{enumerate}[(1)]
\item $\alpha(N)=\alpha'(N)$ iff $\beta(G)$ and $\beta'(G)$ are conjugate by an element of $\Aut(N)$,
\item $\alpha(N)$ is normalized by
$\lambda(G)$ if and only if $\beta(G)$ is contained in the
holomorph $\Hol(N)$ of $N$.
\end{enumerate}

\end{theorem}

As a corollary to the preceding theorem Byott \cite{B}, Proposition 1, obtains the following formula to count Hopf Galois structures on a given finite separable field extension.

\begin{corollary}\label{cor} Let $L/K$ be a separable field extension of degree $g$, $\wL$ its Galois closure, $G=\Gal(\wL/K), G'=\Gal(\wL/L)$. Let $N$ be an abstract group of order $g$ and let $\Hol(N)$ denote the holomorph of $N$. The number $a(N,L/K)$ of Hopf Galois structures of type $N$ on $L/K$ is given by the following formula

$$a(N,L/K)= \dfrac{|\Aut(G,G')|}{|\Aut(N)|} \, b(N,G,G')$$

\noindent where $\Aut(G,G')$ denotes the group of automorphisms of $G$ taking $G'$ to $G'$, $\Aut(N)$ denotes the group of automorphisms of $N$ and $b(N,G,G')$ denotes the number of subgroups $G^*$ of $\Hol(N)$ such that there is an isomorphism from $G$ to $G^*$ taking $G'$ to the stabilizer in $G^*$ of $1_N$.

\end{corollary}

\noindent
{\bf Notation.} In the sequel, $L/K$ will denote a finite separable field extension, $\wL$ the normal closure of $L/K$, $G$ the Galois group $\Gal(\wL/K)$, $G'$ the Galois group $\Gal(\wL/L)$.

\vspace{0.2cm}

In this paper we study Hopf Galois structures on a separable field extension of degree $p^n$, for $p$ an odd prime, $n\geq 2$. In section \ref{n}, we consider the general case $n\geq 2$. In Proposition \ref{pn}, we prove that if a separable field extension of degree $p^n$ has a Hopf Galois structure of cyclic type, then it has no structure of noncyclic type, generalizing Proposition 7 in \cite{CrS}, where the result was established for $n=2$. We note that Kohl proves in \cite{Ko}, Theorem 3.7, that for an odd prime $p$, a cyclic Galois extension of degree $p^n$ admits $p^{n-1}$ Hopf Galois structures, all of cyclic type. In Proposition \ref{pn2}, we generalize his result by determining the number of Hopf Galois structures of cyclic type on a separable field extension of degree $p^n$ such that the Galois group $G=\Gal(\wL/K)$ is a semidirect product $C_{p^n}\rtimes C_D$, with $D$ dividing $p-1$.

\vspace{0.2cm}
In section \ref{3} we study in more detail the case $n=3$. In Proposition \ref{p3} we characterize the transitive groups $G$ of degree $p^3$ for which a separable field extension $L/K$ such that $\Gal(\wL/K)\simeq G$ has some cyclic Hopf Galois structure and determine the number of those. In \cite{CrS}, Theorem 9, we gave the analogous results for degree $p^2$. In subsection \ref{noncy} we consider noncyclic Hopf Galois structures. To this end, we determine the holomorph of each of the four noncyclic groups of order $p^3$. In Theorem \ref{nonab} we prove that if a separable field extension $L/K$ of degree $p^3$ has a nonabelian Hopf Galois structure of type $N$, then it has an abelian Hopf Galois structure whose type has the same exponent as $N$. A related result was obtained in \cite{BC}, theorem 20, where the authors prove that if $L/K$ is a Galois extension with Galois group $G$, a noncyclic abelian group of order $p^n$, $n\geq 3$, then $L/K$ admits a nonabelian Hopf Galois structure.

\vspace{0.2cm}
In Theorem \ref{abin}, we prove that, for $p>3$, the two abelian noncyclic Hopf Galois structures may not occur on the same separable extension of degree $p^3$. For a Galois extension of degree $p^3$, with $p>3$, this result is obtained by applying \cite{C-C-F}, Theorem 1, where the authors prove that if $(N, +)$ is a finite abelian
$p$-group of $p$-rank $m$ where $m + 1 < p$, then every regular abelian subgroup
of the holomorph of $N$ is isomorphic to $N$. The tables in the appendix give the number of Hopf Galois structures of each possible type on a separable extension of degree $27$, for $G$ each of the first 50 transitive groups of degree 27 and illustrates that for $p=3$, all four noncyclic Hopf Galois structures may occur on the same extension. For the Galois case, an example of a regular abelian subgroup of $\Hol({C_3}^3)$ of exponent 9 is given in \cite{C-C-F}, which implies that a Galois extension with Galois group $C_9\times C_3$ has a Hopf Galois structure of type ${C_3}^3$. In our table we may observe that any noncyclic Galois extension of degree $p^3$ has Hopf Galois structures of all four noncyclic types. The fourth row of the table, corresponding to a Galois extension with Galois group the elementary abelian group of order $27$, may be compared with \cite{Ch3}, theorem 7.2, where Childs determines the number of Hopf Galois structures of type $G$, the elementary abelian group of order $p^3$, with $p$ prime $>3$, on a Galois extension with Galois group $G$. Finally, putting together all previous results, we list in Corollary \ref{corf} all possible sets of Hopf Galois structure types on a separable extension of degree $p^3$, for $p>3$ a prime. Corollary \ref{corff} determines the Hopf Galois structure types on a noncyclic abelian extension of degree $p^3$, with $p>3$.

\vspace{0.2cm}
In \cite{B}, Theorem 1, Byott proves that if $G$ is a nilpotent group of order $n$, for each nilpotent group $N$ of order $n$, we obtain $a(N,G)=\prod_{p\mid n} a(\Syl_p(N),\Syl_p(G))$, where $a(N,G)$ (resp. $a(\Syl_p(N),\Syl_p(G))$) denotes the number of Hopf Galois structures of type $N$ (resp. $\Syl_p(N)$) on a Galois extension with Galois group $G$ (resp. $\Syl_p(G)$) and $\Syl_p(G)$ (resp. $\Syl_p(N)$) denotes the $p$-Sylow subgroup of $G$ (resp. of $N$). Hence determining the number of nilpotent Hopf Galois structures on Galois extensions with a nilpotent Galois group reduces to counting Hopf Galois structures on Galois extensions of prime power degree.
In \cite{Z}, Zenouz considers Hopf Galois structures on Galois extensions of order $p^3$ in relation with skew braces (see also \cite{Z2}).

\section{Extensions of degree $p^n$}\label{n}

\begin{proposition}\label{pn}
Let $L/K$ be a separable field extension of degree $p^n$, $p$ an odd prime, $n\geq 2$, $\wL$ its normal closure and $G\simeq \Gal(\widetilde{L}/K)$.
If $L/K$ has a Hopf Galois structure of type $C_{p^n}$, then it has no structure of noncyclic type.
\end{proposition}

\noindent {\it Proof.}  By theorem \ref{theoB}, if $L/K$ has a Hopf Galois structure of type $C_{p^n}$, then $G$ is a transitive subgroup of $\Hol(C_{p^n})$. We shall see that all transitive subgroups of $\Hol(C_{p^n})$ contain an element of order $p^n$. Let us write $\Hol(C_{p^n})$ as $\Z/p^n \Z \rtimes (\Z/p^n \Z)^*$ and let $\sigma$ be a generator of $(\Z/p^n \Z)^*$. The immersion of  $\Hol(C_{p^n})$ in the symmetric group $S_{p^n}$ is given by sending the generator $1$ of $\Z/p^n \Z$ to the $p^n$-cycle $(1,2,\dots,p^n)$ and $\sigma$ to itself, considered as a permutation. We may write

$$\Hol(C_{p^n})=\{ (m,\sigma^j):\, 1\leq m \leq p^n, 1\leq j \leq p^n-p^{n-1} \}.$$

\noindent Let us denote $k:=\sigma(1)$. Then, for $m \in \Z/p^n \Z$, we have $\sigma(m)=mk, \sigma^j(m)=mk^j$. Let $H:=\Hol(C_{p^n})$. Its $p$-Sylow subgroup $\Syl(H)$ has order $p^{2n-1}$ and

$$\Syl(H)=\{ (m,\sigma^{l(p-1)}):\, 1\leq m \leq p^n, 1\leq l \leq p^{n-1} \}.$$

A subgroup $H'$ of $H$ is transitive if and only if $[H':\Stab_H(0) \cap H']=p^n$. Hence if $H'$ is a transitive subgroup of $H$, we have $p^n \mid |H'|$ and $|H'| \mid p^{2n-1}(p-1) = |H|$, therefore $|H'|=p^l d$, with $n\leq l \leq 2n-1$ and $d\mid p-1$. Then $H'$ has a unique $p$-Sylow subgroup $\Syl(H')$ of order $p^l$. We have the following equalities between indices.

$$\begin{array}{ll} & [H':\Stab_H(0) \cap \Syl(H')]=[H':\Syl(H')][\Syl(H'):\Stab_H(0)\cap \Syl(H')], \\
& [H':\Stab_H(0) \cap \Syl(H')]=[H':\Stab_H(0)\cap H'][\Stab_H(0)\cap H':\Stab_H(0)\cap \Syl(H')].\end{array}$$

\noindent Taking into account that $[H':\Stab_H(0)\cap H'] \leq p^n$,  $[\Syl(H'):\Stab_H(0)\cap \Syl(H')]$ is a $p$-power and that $[H':\Syl(H')]$ and $[\Stab_H(0)\cap H':\Stab_H(0)\cap \Syl(H')]$ are divisors of $p-1$, we obtain that $[H':\Stab_H(0) \cap H']=p^n$ if and only if $[\Syl(H'):\Stab_H(0) \cap \Syl(H')]=p^n$. Let us compute $\Stab_H(0) \cap \Syl(H)$.  Since $\Stab_H(0)$ is cyclic of order $p^{n-1}(p-1)$ and generated by $\sigma$, we have $\Stab_H(0) \cap \Syl(H)=\langle \sigma^{p-1}\rangle$, cyclic of order $p^{n-1}$. Now $\Stab_H(0) \cap \Syl(H')=\langle \sigma^{p^a(p-1)}\rangle$, for some integer $a$, $0\leq a\leq n-1$, hence has order $p^{n-1-a}$. Then $|\Syl(H')|=p^{2n-1-a}$.
Writing $\tau:=\sigma^{p-1}$, we have that the order of $(m,\tau^j)$ in $\Syl(H)$ is the maximum of the orders of $m$ in $\Z/p^n\Z$ and $\tau^j$ in $(\Z/p^n \Z)^*$. In particular $(m,\tau^j)$ has order $p^n$ if and only if $m$ has order $p^n$ if and only if $p\nmid m$. Now we observe that the elements in $\Syl(H)$ of order $<p^n$ form a subgroup of $\Syl(H)$, which we denote by $F$. Since $\Syl(H)$ contains $(p^n-p^{n-1})p^{n-1}$ elements of order $p^n$, we have $|F|=p^{2n-2}$. Now $F\cap \Stab_H(0)=F \cap \Syl(\Stab_H(0)) \Rightarrow [F:F \cap \Syl(\Stab_H(0))]=p^{2n-2}/p^{n-1}=p^{n-1}$. We have then obtained that $F$ is not transitive.

Now, if $H'$ has no elements of order $p^n$, we have $\Syl(H') \subset F$, hence $\Syl(H')$ not transitive. This imply $[\Syl(H'):\Syl(H')\cap \Stab_H(0)]<p^n$, hence $H'$ is not transitive. We have then proved that a transitive subgroup of $H$ must contain an element of order~$p^n$.

By \cite{Ko}, Theorem 4.4, if $N$ is any noncyclic group of order $p^n$, then $Hol(N)$ has no elements of order $p^n$. This finishes the proof of the proposition.
 $\Box$

\begin{proposition}\label{pn2}
Let $L/K$ be a separable field extension of degree $p^n$, $p$ an odd prime, $n\geq 2$, $\wL$ its normal closure and $G\simeq \Gal(\widetilde{L}/K)$. Then the number of Hopf Galois structures of cyclic type for $L/K$ is

\begin{enumerate}[1)]
\item $p^{n-1}$ if $L/K$ is a Galois extension with cyclic Galois group;
\item 1 if $G$ is isomorphic to the semidirect product $C_{p^n} \rtimes C_D$ of the cyclic group $C_{p^n}$ of order $p^n$ and a cyclic group $C_D$ of order $D$, with $D$ a divisor of $p-1$.
\end{enumerate}

\end{proposition}

\noindent {\it Proof.} Since an element of order $p^n$ in the symmetric group $S_{p^n}$ is a $p^n$-cycle, we have that a subgroup of $H:=\Hol(C_{p^n})$ containing an element of order $p^n$ is transitive. Taking into account the proof of Proposition \ref{pn}, we obtain that the transitive subgroups of $H$ are exactly those containing an element of order $p^n$.

\begin{enumerate}[1)]
\item In the proof of Proposition \ref{pn}, we saw as well that an element of order $p^n$ of $H$ is of the form $(m,\tau^j$), with $p\nmid m$ and where $\tau$ is an element of order $p^{n-1}$ in $\Aut C_{p^n}$. Since there are $p^{2n-2}(p-1)$ elements of order $p^n$ in $H$ and each cyclic subgroup of order $p^n$ contain $p^{n-1}(p-1)$ such elements, $H$ has $p^{n-1}$ cyclic subgroups of order $p^n$. We obtain then, using Corollary \ref{cor}, that for $G=C_{p^n}$, the number of Hopf Galois structures of type $C_{p^n}$ for $L/K$ is $p^{n-1}$.

\item Let now $H'$ denote a transitive subgroup of $H$ of order $p^n D, D>1$. Then $H'\cap \Stab_H(0)$ has order $D$ and is then a subgroup of order $D$ of $\Stab_H(0)=\langle (0,\sigma) \rangle$. We have then $H'\cap \Stab_H(0)=\langle (0,\sigma^l)\rangle$, with $l=p^{n-1}(p-1)/D$. The subgroups of order $p^n$ of $H$ can be defined as $\langle (1,\tau^j) \rangle$, with $0\leq j \leq p^{n-1}-1$. We have then that the transitive subgroups of $H$ are a product $\langle (1,\tau^j)\rangle \langle (0,\sigma^l) \rangle$.

If $D\mid p-1$, then $\langle (1,\tau^j) \rangle$ is the $p$-Sylow subgroup of $H'$, hence it is normal in $H'$. Now $(0,\sigma^l)(1,\tau^j)(0,\sigma^{-l})=(k^l,\tau^j)$, where $k=\sigma(1)$. If $\tau^j \neq \Id$, then $(k^l,\tau^j) \in \langle (1,\tau^j) \rangle \Leftrightarrow (k^l,\tau^j)=(1,\tau^j)^N$, with $N$ such that $jN \equiv j \pmod{p^{n-1}}$. But then the first component of $(1,\tau^j)^N$ is

$$1+k^{j(p-1)}+ \dots +k^{j(p-1)(N-1)} =\dfrac{1-k^{(p-1)jN}}{1-k^{(p-1)j}} \equiv 1 \pmod{p^n}$$

\noindent and it cannot be congruent to $k^l$ modulo $p^n$. We have then that the only transitive subgroup of $H$ of order $p^n D$ is  $H'=\langle (1,\Id),(0,\sigma^l) \rangle$. Now since $\langle (1,\Id)\rangle$ is the unique $p$-Sylow subgroup of $H'$, an automorphism of $H'$ must send $(1,\Id)$ to $(m,\Id)$, for some $m$ not divisible by $p$. Since $(0,\sigma^l)(m,\Id)(0,\sigma^{-l})=(mk^l,\Id)$, an automorphism of $H'$ sending $\langle (0,\sigma^l) \rangle$ to itself must send the element $(0,\sigma^l)$ to itself. We obtain then, using Corollary \ref{cor}, that for $G=C_{p^n}\rtimes C_D$, with $D\mid p-1$, the number of Hopf Galois structures of type $C_{p^n}$ for $L/K$ is $1$. $\Box$
\end{enumerate}

\section{Extensions of degree $p^3$}\label{3}
\subsection{Hopf Galois structures of cyclic type}

We determine now exactly the number of Hopf Galois structures of cyclic type for separable field extensions of degree $p^3$, $p$ an odd prime.

\begin{proposition}\label{p3}
Let $L/K$ be a separable field extension of degree $p^3$, $p$ an odd prime, $\wL$ its normal closure and $G\simeq \Gal(\widetilde{L}/K)$. Then $L/K$ has a Hopf Galois structure of type $C_{p^3}$ if and only if $G$ is isomorphic to the semidirect product $C_{p^3} \rtimes C_D$ of the cyclic group $C_{p^3}$ of order $p^3$ and a cyclic group $C_D$ of order $D=p^a d$, with $0\leq a \leq 2, d$ a divisor of $p-1$.

Moreover the number of Hopf Galois structures of type $C_{p^3}$ for $L/K$ is 1 if $d>1$, $p^2$ if $D=1$ and $p^{3-a}$ if $D=p^a$, with $a=1$ or $2$.

\end{proposition}

\noindent {\it Proof.} In the proof of Proposition \ref{pn2}, we have seen that the transitive subgroups of $H:=\Hol(C_{p^n})$ are a product $\langle (1,\tau^j)\rangle \langle (0,\sigma^l) \rangle$. We have to see when this product is a subgroup $H'$ of $H$.

\begin{enumerate}[1)]
\item If $D\mid p-1$, we have seen in Proposition \ref{pn2} that the only transitive subgroup of $H$ of order $p^3 D$ is  $\langle (1,\Id),(0,\sigma^l) \rangle \simeq C_{p^3} \rtimes C_D$  and that for $G=C_{p^3}\rtimes C_D$, with $D\mid p-1$, there is exactly one Hopf Galois structure of type $C_{p^3}$ for $L/K$.
\item If $D=p^a$, with $a=1$ or $2$, then $\sigma^l=\tau^{p^{2-a}}$ and the product  $\langle (1,\tau^j)\rangle \langle (0,\tau^{p^{2-a}}) \rangle$ contains $p^a$ cyclic subgroups of order $p^3$, namely $\langle (1,\tau^{j+\lambda p^{2-a}})\rangle$, with $0\leq \lambda < p^a$. Therefore we have $p^{2-a}$ possible transitive subgroups of $H$ of order $p^{3+a}$, namely the products
    $\langle (1,\tau^j)\rangle \langle (0,\tau^{p^{2-a}})\rangle$, with $0\leq j < p^{2-a}$.

\begin{enumerate}[i)]
\item
    If $a=2$, there is just one product $\langle (1,\Id)\rangle \langle (0,\tau)\rangle$ which is a group isomorphic to $C_{p^3} \rtimes C_{p^2}$ (where the morphism from $C_{p^2}$ to $\Aut(C_{p^3})$ giving the action of $C_{p^2}$ on $C_{p^3}$ is into). Since $(0,\tau)(1,\Id)(0,\tau^{-1})=(\tau(1),\Id)=(1,\Id)^{\tau(1)}$, an automorphism of $\langle (1,\Id), (0,\tau)\rangle$ sending $\langle (0,\tau)\rangle$ to itself must send the element $(0,\tau)$ to itself. Now $(0,\tau)(m,\tau^j)(0,\tau^{-1})=(m\tau(1),\tau^j)$ which equals $(m,\tau_j)^{\tau(1)}$ if and only if $p\mid j$. We have then that an automorphism of $\langle (1,\Id), (0,\tau)\rangle$ sending $\langle (0,\tau)\rangle$ to itself is given by

    $$(1,\Id) \mapsto (m,\tau^j), \text{\ with \  } p\nmid m, p\mid j, \quad (0,\tau) \mapsto (0,\tau)$$

    \noindent and hence there are $(p^3-p^2)p$ such automorphisms. The number of structures of cyclic type for $L/K$ when $G\simeq C_{p^3} \rtimes C_{p^2}$ is then $p$.
\item
    If $a=1$, there are $p$ products  $\langle (1,\tau^j)\rangle \langle (0,\tau^{p})\rangle$, with $0\leq j < p$. We may assume $\tau(1)=p+1$ and obtain $\sum_{c=0}^{p^2} \tau^{cj}(1) \equiv 1+p^2 \pmod{p^3}$ which implies

    \begin{equation}\label{eq1}
    (0,\tau^{p})(1,\tau^j)(0,\tau^{-p})=(\tau^{p}(1),\tau^j)=(1,\tau^j)^{\tau^p(1)}.
    \end{equation}

\noindent Hence the products considered are subgroups of $H$ isomorphic to $C_{p^3} \rtimes C_p$. Taking into account (\ref{eq1}), we obtain
that an automorphism of $\langle (1,\tau^j), (0,\tau^p)\rangle$ sending $\langle (0,\tau^p)\rangle$ to itself must send the element $(0,\tau^p)$ to itself and that such automorphisms are given by

    $$(1,\tau^j) \mapsto (1,\tau^j)^m(0,\tau^{\lambda p}), \text{\ with \  } p\nmid m, 0\leq \lambda <p, \quad (0,\tau) \mapsto (0,\tau).$$

\noindent The number of these automorphisms is $(p^3-p^2)p$ and hence the number of structures of cyclic type for $L/K$ when $G\simeq C_{p^3} \rtimes C_{p}$ is $p^2$.

\end{enumerate}

\item If $D$ and $p-1$ are not coprime, we may write $D=p^a d$, with $d>1$, $d\mid p-1$. If $G$ is a transitive subgroup of order $p^{3+a} d$, then it has a unique $p$-Sylow subgroup which is a transitive subgroup of $H$ of order $p^{3+a}$.
\begin{enumerate}[i)]
\item
    If $a=2$, this $p$-Sylow subgroup must be $\langle (1,\Id)\rangle \langle (0,\tau)\rangle$ and then \linebreak $G=\langle (1,\Id)\rangle \langle (0,\sigma^l)\rangle$, with $l=(p-1)/d$. Since $(0,\sigma^l)(1,\Id)(0,\sigma^{-l})=(\sigma^l(1),\Id)=(1,\Id)^{\sigma^l(1)}$, $G$ is indeed a group, isomorphic to $C_{p^3}\rtimes C_D$. As in the preceding cases, an automorphism of $G$ sending $\langle(0,\sigma^l)\rangle$ to itself must send the element $(0,\sigma^l)$ to itself. If there were such an automorphism of $G$ sending $(1,\Id)$ to $(1,\tau^j)$, we would have $(0,\sigma^l)(1,\tau^j)(0,\sigma^{-l})=(\sigma^l(1),\tau^j)=(1,\tau^j)^{\sigma^l(1)}$ but we have seen in the proof of 2) of Proposition \ref{pn2} that $(\sigma^l(1),\tau^j)\in \langle(1,\tau^j)\rangle$ implies $j=0$. We have then that the automorphisms of $G$ sending $\langle(0,\sigma^l)\rangle$ to itself are given by $(1,\Id) \mapsto (m,\Id)$, with $p\nmid m$, $(0,\sigma^l)\mapsto (0,\sigma^l).$ We have then $p^3-p^2$ such automorphisms and there is exactly one structure of cyclic type for $L/K$ when $G\simeq C_{p^3} \rtimes C_{D}$, with $D=p^2d, d\mid p-1, d>1$.
\item
    If $a=1$, this $p$-Sylow subgroup must be $\langle (1,\tau^j),(0,\tau^{p})\rangle$, for some $j$ with $0\leq j < p$ and is normal in $G$. Now, if $j\neq 0$, $(0,\sigma^l)(1,\tau^j)(0,\sigma^{-l})=(\sigma^l(1),\tau^j)=(1,\tau^j)^N(0,\tau^p)^M \Rightarrow Nj+Mp \equiv j \pmod{p^2} \Rightarrow N \equiv 1 \pmod{p}$. As in the proof of 2) in Proposition \ref{pn2}, the first component of $(1,\tau^j)^N(0,\tau^p)^M $ would then be congruent to 1, modulo $p$ and cannot equal $\sigma^l(1)$. We have then $G=\langle (1,\Id)\rangle \langle (0,\sigma^l)\rangle$, with $l=p(p-1)/d$. As in the case $a=2$, there are $p^3-p^2$ automorphisms of $G$ sending $\langle(0,\sigma^l)\rangle$ to itself and there is exactly one structure of cyclic type for $L/K$ when $G\simeq C_{p^3} \rtimes C_{D}$, with $D=pd, d\mid p-1, d>1$. $\Box$
\end{enumerate}
\end{enumerate}

\subsection{Hopf Galois structures of noncyclic type}\label{noncy}

Let $p$ be an odd prime. There are exactly 5 groups of order $p^3$ up to isomorphism: three abelian ones $C_{p^3}, C_{p^2}\times C_p, C_p\times C_p\times C_p$ and two non abelian ones, the Heisenberg group $H_p$ and a group $G_p$ of exponent $p^2$, defined as follows.

$$H_p:= \left\{ \left( \begin{array}{ccc} 1&a&b\\0&1&c\\ 0&0&1 \end{array} \right) : a,b,c \in \F_p \right\} \subset \GL(3,\F_p),$$

$$G_p:= \left\{ \left(\begin{array}{cc} 1+pb & a \\ 0 & 1 \end{array} \right) : a, b \in \Z/p^2 \Z \right\} \subset \GL(2,\Z/p^2 \Z),$$

\noindent where $b$ is taken modulo $p$ (see \cite{C}). Let us note that all nontrivial elements in $H_p$ have order $p$ while $G_p$ has $p^3-p^2$ elements of order $p^2$, those with $a\not \equiv 0 \pmod{p}$.

We shall study the Hopf Galois structures of noncyclic type for separable extensions of degree $p^3$, with $p$ an odd prime, using Theorem \ref{theoB}. Since for a group $N$, we have $\Hol(N)=N\rtimes \Aut N$, we determine first the automorphism group of the four noncyclic groups of order $p^3$.

\begin{enumerate}[1)]
\item Let us write $C_{p^2} \times C_p=\langle a \rangle \times \langle b \rangle$. The elements of order $p^2$ in $C_{p^2} \times C_p$ are those of the form $a^i b^j$, with $p\nmid i$, there are $(p^2-p)p$ such elements. The elements of order $p$ in $C_{p^2} \times C_p$ are those of the form $a^{pk} b^l$, with $p$ not dividing both $k$ and $l$, there are $p^2-1$ such elements. If $\varphi$ is an automorphism of $C_{p^2} \times C_p$, we have $\varphi(a)=a^i b^j$, with $p\nmid i$. Now, since $\varphi(b) \not \in \langle\varphi(a)\rangle$, we have $\varphi(b)=a^{pk} b^l$, with $p \nmid l$. We obtain then $|\Aut(C_{p^2} \times C_p)|=p^3(p-1)^2$.
\item $\Aut(C_p^3) \simeq \GL(3,\F_p)$ and $|\GL(3,\F_p)|=(p^3-1)(p^3-p)(p^3-p^2)=p^3(p-1)^3(p+1)(p^2+p+1)$
\item $H_p$ is generated by $A=\left(\begin{smallmatrix} 1&1&0\\0&1&0\\0&0&1 \end{smallmatrix}\right)$ and $C=\left(\begin{smallmatrix} 1&0&0\\0&1&1\\0&0&1 \end{smallmatrix}\right)$ satisfying $AC=BCA$, where $B=\left(\begin{smallmatrix} 1&0&1\\0&1&0\\0&0&1 \end{smallmatrix}\right)$. Moreover the centre of $H_p$ is equal to $\langle B \rangle$ and $H_p/\langle B \rangle \simeq C_p\times C_p$.
     This gives an epimorphism $\Aut(H_p) \rightarrow \Aut(C_p\times C_p) \simeq \GL(2,\F_p)$. An automorphism $\varphi$ of $H_p$ is then given by $\varphi(A)=B^iA^rC^s, \varphi(C)=B^jA^tC^v$, with $p \nmid rv-st$ and we get $\varphi(B)=B^{rv-st}$. We have then $|\Aut(H_p)|=p^2|\GL(2,\F_p)|=p^2(p^2-1)(p^2-p)=p^3(p-1)^2(p+1).$
\item $G_p$ is generated by $M=\left(\begin{smallmatrix} 1&1\\0&1 \end{smallmatrix}\right)$ and $N=\left(\begin{smallmatrix} 1+p&0\\0&1 \end{smallmatrix}\right)$, where $M$ has order $p^2$, $N$ has order $p$ and $NM=M^{p+1}N$. An automorphism $\varphi$ of $G_p$ is given by $\varphi(M)=\left(\begin{smallmatrix} 1+pb&a\\0&1 \end{smallmatrix}\right)$, with $a\not \equiv 0 \pmod{p}$, and $\varphi(N)=\left(\begin{smallmatrix} 1+p&pc\\0&1 \end{smallmatrix}\right)$. Hence $|\Aut(G_p)|=p^3(p-1)$. For further use we note that $G_p$ is also generated by $M$ and $N_j=\left(\begin{smallmatrix} 1+jp&0\\0&1 \end{smallmatrix}\right)$, with $1\leq j \leq p-1$ satisfying $N_jM=M^{jp+1}N_j$.
\end{enumerate}

The following theorem relates Hopf Galois structures of abelian and nonabelian types.

\begin{theorem}\label{nonab} Let $p$ be an odd prime. For a separable field extension  $L/K$ of degree $p^3$, the following implications hold.

\begin{enumerate}[1)]
\item $L/K$ has a Hopf Galois structure of type $H_p \Rightarrow L/K$ has a Hopf Galois structure of type $C_p^3$,
\item $L/K$ has a Hopf Galois structure of type $G_p \Rightarrow L/K$ has a Hopf Galois structure of type $C_{p^2} \times C_p$.
\end{enumerate}
\end{theorem}

\noindent {\it Proof.} 1) We identify $C_p$ with the additive group of the field $\F_p$ and the symmetric group $S_{p^3}$ with the group of permutations of $\F_p^3$. Then the action of $C_p^3$ on itself by left translation induces a monomorphism

$$\begin{array}{lcll} \lambda_{C_p^3}:&C_p^3 &\rightarrow & \Perm(\F_p^3) \\
 & (l,m,n) & \mapsto & t_{(l,m,n)} : (a,b,c)\mapsto (a+l,b+m,c+n) \end{array}.$$

\noindent We consider the bijection

$$b:H_p \rightarrow \F_p^3, \left( \begin{smallmatrix} 1&a&b\\0&1&c\\ 0&0&1 \end{smallmatrix} \right) \mapsto (a,b-ac/2,c)$$

\noindent
and the isomorphism $\beta:\Perm(H_p) \rightarrow \Perm(\F_p^3)$ given by $\beta(\sigma)=b\sigma b^{-1}$. The action of $H_p$ on itself by left translation induces a monomorphism

$$\lambda_{H_p}:H_p \rightarrow \Perm(\F_p^3)$$

\noindent given by

$$\begin{array}{l} \lambda_{H_p}(A)(a,b,c)=(a+1,b+c/2,c), \quad \lambda_{H_p}(B)(a,b,c)=(a,b+1,c), \\ \lambda_{H_p}(C)(a,b,c)=(a,b-a/2,c+1).\end{array}$$

\noindent We want to prove the inclusion $\Norm_{\Perm(\F_p^3)}(\lambda_{H_p}(H_p)) \subset \Norm_{\Perm(\F_p^3)}(\lambda_{C_p^3}(C_p^3))$. By Theorem \ref{theoB}, this implies 1). We observe that $\lambda_{H_p}(B)=\lambda_{C_p^3}(0,1,0)$. Note that, since $C_p^3$ is abelian, $\lambda_{C_p^3}$ coincides with the monomorphism $\rho_{C_p^3}$ induced by the action of $C_p^3$ on itself by right translation. Hence $\lambda_{H_p}(B)=\rho_{C_p^3}(0,1,0) \in \Norm_{\Perm(\F_p^3)}(\lambda_{C_p^3}(C_p^3))$. Moreover

$$\lambda_{H_p}(A) t_{(l,m,n)}\lambda_{H_p}(A)^{-1}=t_{(l,m+n/2,n)},\quad
\lambda_{H_p}(C)t_{(l,m,n)}\lambda_{H_p}(C)^{-1}=t_{(l,m-l/2,n)}.$$

\noindent hence $\lambda_{H_p}(A)$ and $\lambda_{H_p}(C)$ normalize $\lambda_{C_p^3}(C_p^3)$. We have then

$$\lambda_{H_p}(H_p)\subset \Norm_{\Perm(\F_p^3)}(\lambda_{C_p^3}(C_p^3)).$$

We want to see now that the image of $\Aut(H_p)$ by the inclusion in $\Perm(H_p)$ followed by the isomorphism $\beta: \Perm(H_p) \rightarrow \Perm(\F_p^3)$ is contained in $\Norm_{\Perm(\F_p^3)}(\lambda_{C_p^3}(C_p^3))$. It is enough to consider automorphisms of the two following forms

$$\begin{array}{cccc} \varphi_1:& A & \mapsto & B^i A \\ & C & \mapsto & B^j C \end{array} \quad , \quad \begin{array}{cccc} \varphi_2:& A & \mapsto &  A^rC^s \\ & C & \mapsto & A^t C^v \end{array}, \text{with \ } p \nmid rv-st.$$

\noindent Indeed, for $i,j$ running over all integers modulo $p$, $\varphi_1$ runs over all elements in the kernel of the epimorphism $\Aut(H_p) \rightarrow \GL(2,\F_p)$ and for $r,s,t,v$ running over all integers modulo $p$, with $p \nmid rv-st$, $\varphi_2$ runs over all cosets of this kernel in $\Aut(H_p)$.

We denote again by $\varphi_1$ and $\varphi_2$ their images in $\Perm(\F_p^3)$. Note that $\left(\begin{smallmatrix} 1&a&b\\0&1&c \\ 0&0&1 \end{smallmatrix}\right)=C^cA^aB^b$. As an element in $\Aut(H_p)$, $\varphi_1$ satisfies $\varphi_1\left(\begin{smallmatrix} 1&a&b\\0&1&c \\ 0&0&1 \end{smallmatrix}\right)=\left(\begin{smallmatrix} 1&a&b+ia+jc\\0&1&c \\ 0&0&1\end{smallmatrix}\right)$. Taking into account the isomorphism $\beta$ between $\Perm(H_p)$ and $\Perm(\F_p^3)$, we obtain $\varphi_1((a,b,c))=(a,b+ia+jc,c)$ and $\varphi_1 t_{(l,m,n)} \varphi_1^{-1}=t_{(l,m+il+jn,n)}$. As an element in $\Aut(H_p)$, $\varphi_2$ satisfies $\varphi_2\left(\begin{smallmatrix} 1&a&b\\0&1&c \\ 0&0&1 \end{smallmatrix}\right)=\left(\begin{smallmatrix} 1&a'&b'\\0&1&c' \\ 0&0&1\end{smallmatrix}\right)$ with

$$\begin{array}{l} a'=ct+ar \\
b'=-\dfrac{(c-1)ctv+(a-1)ars} 2 +brv-bst+c^2tv+acst+a^2rs \\
c'=cv+as \end{array}.$$

\noindent Taking into account the isomorphism $\beta$ between $\Perm(H_p)$ and $\Perm(\F_p^3)$, we obtain $\varphi_2((a,b,c))=(ct+ar,\dfrac 1 2 ctv +\dfrac 1 2 ars+brv-bst,cv+as)$ and $\varphi_2 t_{(l,m,n)} \varphi_2^{-1}=t_{(l',m',n')}$ with

$$\begin{array}{l} l'=nt+lr \\
m'=\dfrac 1 2 (ntv+lrs)+mrv-mst \\
n'=nv+ls \end{array}.$$

We have then

$$\beta(\Aut(H_p)) \subset \Norm_{\Perm(\F_p^3)}(\lambda_{C_p^3}(C_p^3)).$$

\vspace{0.3cm}

\noindent 2) Let $C_{p^2}\times C_p=<a,b>$. We consider the following automorphisms of $C_{p^2}\times C_p$.

$$\begin{array}{l} \begin{array}{cccl} \varphi_1:& a &\mapsto &a^{p+1} \\ & b& \mapsto & b \end{array}\quad \quad
\begin{array}{cccl} \varphi_2:& a &\mapsto &a \\ & b& \mapsto & a^p b \end{array} \quad \quad
\begin{array}{cccl} \varphi_3:& a &\mapsto &ab \\ & b& \mapsto & b \end{array} \\
\begin{array}{cccl} \psi_1:& a &\mapsto &a^i \\ & b& \mapsto & b \end{array} \quad \quad \quad
\begin{array}{cccl}  \psi_2:& a &\mapsto &a \\ & b& \mapsto & b^l\end{array}\end{array}
$$

\noindent where $i$ has order $p-1$ modulo $p^2$ and $l$ has order $p-1$ modulo $p$. The automorphisms $\varphi_1, \varphi_2, \varphi_3$  have order $p$ and satisfy $\varphi_1 \varphi_2=\varphi_2 \varphi_1, \varphi_1\varphi_3=\varphi_3\varphi_1, \varphi_2\varphi_3=\varphi_1\varphi_3\varphi_2$. We have then that the subgroup $\langle \varphi_1, \varphi_2, \varphi_3\rangle$ of $\Aut(C_{p^2} \times C_p)$ is isomorphic to $H_p$ and is then a $p$-Sylow subgroup of $\Aut(C_{p^2} \times C_p)$. The automorphisms $\psi_1$ and $\psi_2$ have order $p-1$ and commute with each other. Moreover $\psi_1\varphi_1\psi_1^{-1}=\varphi_1, \psi_1\varphi_2\psi_1^{-1}=\varphi_2^i,
\psi_1\varphi_3\psi_1^{-1}=\varphi_3^{i'}, \psi_2\varphi_1\psi_2^{-1}=\varphi_1, \psi_2\varphi_2\psi_2^{-1}=\varphi_2^{l'},
\psi_2\varphi_3\psi_2^{-1}=\varphi_3^{l}$, where $i'$ is the inverse of $i$ modulo $p^2$ and $l'$ is the inverse of $l$ modulo $p$. Hence $\varphi_1, \varphi_2, \varphi_3, \psi_1, \psi_2$ generate $\Aut(C_{p^2} \times C_p)$ and the $p$-Sylow subgroup is unique. We have then that $\Hol(C_{p^2}\times C_p)$ is generated by $a,b,\varphi_1, \varphi_2, \varphi_3, \psi_1, \psi_2$.

We consider now the elements $\mathcal{M}:=(a,\varphi_2), \mathcal{N}:=(b,\varphi_1^{-1})$ of $\Hol(C_{p^2} \times C_p)=(C_{p^2} \times C_p)\rtimes \Aut(C_{p^2} \times C_p)$. We may check that $\mathcal{M}$ has order $p^2$, $\mathcal{N}$ has order $p$ and $\mathcal{N}\mathcal{M}=\mathcal{M}^{-2p+1}N$, hence $\langle \mathcal{M},\mathcal{N} \rangle \simeq G_p$. Moreover, for given $i,j$, $(a,\varphi_2)^{(1-p)i}(b,\varphi_1^{-1})^j$ sends $1_{C_{p^2} \times C_p}$ to $a^ib^j$, hence $\langle \mathcal{M},\mathcal{N} \rangle$ is a transitive subgroup of $\Sym(C_{p^2} \times C_p)$. Now the relations

$$\begin{array}{ll} (1,\varphi_1)(a,\varphi_2)(1,\varphi_1^{-1})= (a,\varphi_2)^{p+1} & (1,\varphi_1)(b,\varphi_1^{-1})(1,\varphi_1^{-1})= (b,\varphi_1^{-1}) \\
(1,\varphi_2)(a,\varphi_2)(1,\varphi_2^{-1})= (a,\varphi_2)& (1,\varphi_2)(b,\varphi_1^{-1})(1,\varphi_2^{-1})= (a,\varphi_2)^p(b,\varphi_1^{-1}) \\
(1,\varphi_3)(a,\varphi_2)(1,\varphi_3^{-1})=(b,\varphi_1^{-1}) (a,\varphi_2)^{p+1}& (1,\varphi_3)(b,\varphi_1^{-1})(1,\varphi_3^{-1})=(b,\varphi_1^{-1}) \\
(1,\psi_1)(a,\varphi_2)(1,\psi_1^{-1})=(a,\varphi_2)^i & (1,\psi_1)(b,\varphi_1^{-1})(1,\psi_1^{-1})= (b,\varphi_1^{-1})\end{array}$$

\noindent with $i$ as in the definition of $\psi_1$, imply that the normalizer of $\langle \mathcal{M},\mathcal{N} \rangle$ in the normalizer of $C_{p^2}\times C_p$ has order at least $p^6(p-1)$, hence the normalizer of $\langle \mathcal{M},\mathcal{N} \rangle$ in $\Sym(C_{p^2}\times C_p)$ is contained in the normalizer of $C_{p^2}\times C_p$ in $\Sym(C_{p^2}\times C_p)$. $\Box$

\vspace{0.5cm}
In the tables in the appendix, we present for the first 50 transitive groups $G$ of degree 27, the number of Hopf Galois structures of each of the five possible types on a separable extension $L/K$ of degree 27 such that $G \simeq \Gal(\wL/K)$ as well as the total number of these Hopf Galois structures. We have made this computation by implementing in Magma the formula in Corollary \ref{cor}. We observe that for several transitive groups $G$, in particular for all noncyclic $G$ of order 27, all four noncyclic types of Hopf Galois structures appear. In the next proposition, we prove that this fact only occurs for $p=3$.

\begin{theorem}\label{abin} Let $p$ be an odd prime, $p>3$. Let $L/K$  be a separable field extension  of degree $p^3$. If $L/K$ has a Hopf Galois structure of type $C_{p^2}\times C_p$, then it has no Hopf Galois structure of type $C_p \times C_p \times C_p$.
\end{theorem}

\noindent
{\it Proof.} Let $\wL$ denote the Galois closure of $L/K$ and $G=\Gal(\wL/K)$. If $L/K$ has some Hopf Galois structure of type $C_{p^2}\times C_p$ and some Hopf Galois structure of type $C_p \times C_p \times C_p$, then $G$ is isomorphic to a transitive subgroup of $\Hol(C_{p^2}\times C_p)$ and also to a transitive subgroup of $\Hol(C_p \times C_p \times C_p)$. Then if $Q$ is a $p$-Sylow subgroup of $G$, then $Q$ is a subgroup of some $p$-Sylow subgroup of $\Hol(C_{p^2}\times C_p)$ and also a  subgroup of some $p$-Sylow subgroup of $\Hol(C_p \times C_p \times C_p)$.
With the notations in the second part of the proof of Theorem \ref{nonab}, $P_1:=\langle a,b\rangle \rtimes \langle\varphi_1,\varphi_2,\varphi_3 \rangle$ is a $p$-Sylow subgroup of
$\Hol(C_{p^2}\times C_p)$ and, taking into account the relations of $\psi_1$ and $\psi_2$, it is normal in $\Hol(C_{p^2}\times C_p)$, hence it is unique. Now, since $p^3$ divides exactly $|\Aut(C_p \times C_p \times C_p)|$, a $p$-Sylow subgroup of $\Hol(C_p \times C_p \times C_p)$ has order $p^6$ and is isomorphic to $P_2:=\F_p^3 \rtimes H_p$. We will check that $P_1$ and $P_2$ are not isomorphic by computing the number of elements of order $p$ in each of them.

Let $(x,\varphi)$, with $x=a^jb^k, \varphi=\varphi_1^{m_1}\varphi_2^{m_2}\varphi_3^{m_3}$, be an element in $P_1$. Then $(x,\varphi)^p=(x\varphi(x) \dots \varphi^{p-1}(x),\Id)$. Now, we have

$$\varphi(a)=a^{1+p(m_1+m_2m_3)} b^{m_3}, \quad \varphi(b)=a^{m_2p} b$$

\noindent and, for $k$ a nonnegative integer,

$$\varphi^{k}(a)=a^{1+p(km_1+\frac{(k+1)k} 2 m_2m_3)} b^{km_3}, \quad \varphi^k(b)=a^{km_2p} b.$$

\noindent We obtain then $b\varphi(b) \dots \varphi^{p-1}(b)=a^{m_2p\frac{p(p-1)} 2} b^p=1$ and, with $S=\sum_{k=1}^{p-1} \frac{k(k+1)} 2$,

$$a\varphi(a) \dots \varphi^{p-1}(a)=a^{p+p(\frac{p(p-1)} a m_1+Sm_2m_3)}b^{\frac {p(p-1)} 2 m_3}=a^{p(1+Sm_2m_3)}.$$

\noindent Now

$$S=\sum_{k=1}^{p-1} \dfrac{k(k+1)} 2= \dfrac 1 2 \dfrac{(p-1)p(2p-1)}6 + \dfrac 1 2 \dfrac{p(p-1)} 2.$$

\noindent If $p=3$, we have $S=4$. If $p\equiv 1 \pmod{3}$, writing $p=3q+1$ and taking into account that $q$ is even, we obtain $S=(q/2) (3q+1)(3q+2)$ a multiple of $p$. If $p\equiv 2 \pmod{3}$, writing $p=3q+2$ and taking into account that $q$ is odd, we obtain $S=((q+1)/2) (3q+1)(3q+2)$ a multiple of $p$. We obtain then that, for $p\neq 3$, $S$ is a multiple of $p$, hence $a\varphi(a) \dots \varphi^{p-1}(a)=a^{p}$. Since $a$ and $b$ commute with each other, we obtain that $(a^jb^k,\varphi)$ has order $p$ if and only if $p\mid j$. Hence there are $p^5-1$ elements of order $p$ in $P_1$. In the case $p\nmid j$, since $(a^p,\Id)$ has order $p$, $(a^jb^k,\varphi)$ has order $p^2$.

We identify $(C_p)^3$ with the additive group of $\F_p^3$ and consider $v=(x,y,z) \in \F_p^3, M=\left( \begin{smallmatrix} 1&a&b\\0&1&c\\ 0&0&1 \end{smallmatrix} \right)\in H_p$. Then $(v,M)^p=(v+Mv+\dots+M^{p-1}v,\Id)$. Now

$$M^k=\left( \begin{array}{ccc} 1&ka& kb+\dfrac{k(k-1)} 2 ac \\ 0 & 1 & kc\\ 0&0&1 \end{array} \right)$$

\noindent and

$$\Id+M+\dots+M^{p-1}=\left( \begin{array}{ccc} p&\dfrac{p(p-1)} 2a& \dfrac{p(p-1)} 2 b+S' ac \\ 0 & p & \dfrac{p(p-1)}2 c\\ 0&0&p \end{array} \right)$$

\noindent where

$$S'=\sum_{k=1}^{p-1} \dfrac{k(k-1)} 2= \dfrac 1 2 \dfrac{(p-1)p(2p-1)}6 - \dfrac 1 2 \dfrac{p(p-1)} 2.$$

\noindent If $p=3$, we have $S'=1$. If $p\equiv 1 \pmod{3}$, writing $p=3q+1$ and taking into account that $q$ is even, we obtain $S'=(q/2) (3q+1)(3q-1)$ a multiple of $p$. If $p\equiv 2 \pmod{3}$, writing $p=3q+2$ and taking into account that $q$ is odd, we obtain $S'=((3q-1)/2) (3q+2)q$ a multiple of $p$. We obtain then that, for $p\neq 3$, $S'$ is a multiple of $p$, hence $\Id+M+\dots+M^{p-1}=0$ and all elements in $P_2$ different from $(0,\Id)$ have order $p$.

We have then that the $p$-Sylow group $Q$ of $G$ has order at most $p^5$ and since all nontrivial elements in $P_2$ have order $p$, the same has to hold for $Q$. Now the subgroup $P_1':=\langle a^p,b,\varphi_1,\varphi_2,\varphi_3 \rangle$ of $P_1$ has order $p^5$ and contains all elements of order $p$ in $P_1$. We have then $Q\subset P_1'$ and hence $G\subset \langle a^p,b\rangle\rtimes \Aut(C_{p^2} \times C_p)$. Now the orbit of $1$ under the action of $\langle a^p,b\rangle$ by translation is $X:=\{ a^{pj} b^k\, : \, 0\leq j,k <p \}$ and the orbit of any element in $X$ by the action of $\Aut(C_{p^2} \times C_p)$ is contained in $X$. Hence $\langle a^p,b\rangle\rtimes \Aut(C_{p^2} \times C_p)$ is not transitive neither is then $G$ and we obtain a contradiction. $\Box$

\begin{corollary}\label{corf}
Let $p>3$ be a prime. Let $L/K$ be a separable field extension of degree $p^3$. Then either $L/K$ has no Hopf Galois structures or the set of types of Hopf Galois structures on $L/K$ is one the following:

$$\{C_{p^3} \}, \quad \{(C_p)^3 \}, \quad \{C_{p^2}\times C_p\},\quad \{(C_p)^3,H_p \}, \quad\{C_{p^2}\times C_p,G_p\}.$$

\end{corollary}

\noindent {\it Proof.} By Proposition \ref{pn} and Theorems \ref{nonab} and \ref{abin} the possibilities not considered in the corollary do not occur. Let us see now that all cases listed do occur by exhibit an example. We denote again by $\wL$ a Galois closure of $L/K$ and $G=\Gal(\wL/K)$. For $G=S_{p^3}$, the whole symmetric group in $p^3$ letters, $L/K$ has no Hopf Galois structures, by \cite{G-P} Corollary 4.8. If $G=C_{p^3}$, then $L/K$ has only Hopf Galois structures of cyclic type by Proposition \ref{pn}. If $G=\Hol((C_p)^3)$, then $L/K$ has only Hopf Galois structures of type $(C_p)^3$ since $|G|>|\Hol(N)|$ for any other group $N$ of order $p^3$. If $G=\Hol(C_{p^2}\times C_p)$, then $G$ cannot be embedded in $\Hol(G_p)$, since $|G|>|\Hol(G_p)|$, and by Proposition \ref{pn} and Theorems \ref{nonab} and \ref{abin}, $L/K$ may not have Hopf Galois extensions of the remaining types. If $G=H_p$, then by Theorem \ref{nonab}, $L/K$ has also Hopf Galois structures of type $(C_p)^3$ and by Proposition \ref{pn} and Theorems \ref{nonab} and \ref{abin}, $L/K$ may not have Hopf Galois extensions of the remaining types. Similarly if $G=G_p$, then $L/K$ has Hopf Galois structures precisely of types $G_p$ and $C_{p^2}\times C_p$. $\Box$

\vspace{0.2cm}
In the next corollary, we determine the types of Hopf Galois structures on an abelian noncyclic extension of degree $p^3$, $p>3$ a prime.

\begin{corollary}\label{corff} Let $p>3$ be a prime. Let $L/K$ be a Galois extension with Galois group~$G$.
\begin{enumerate}[1)]
\item If $G\simeq C_{p^2} \times C_p$, then the Hopf Galois structures on $L/K$ are exactly of types  $C_{p^2} \times C_p$ and $G_p$;
\item If $G\simeq C_{p}^3$, then the Hopf Galois structures on $L/K$ are exactly of types  $C_{p}^3$ and $H_p$.
\end{enumerate}
\end{corollary}

\noindent
{\it Proof.} By Corollary \ref{corf}, $L/K$ has Hopf Galois structures either only of type $G$ or of type $G$ and type the nonabelian group of order $p^3$ with same exponent as $G$. By \cite{BC}, theorem 20, $L/K$ has some nonabelian Hopf Galois structure. Hence the statement is proved.

\newpage
\addtolength{\textwidth}{1cm}
\begin{landscape}
\section*{Appendix: Hopf Galois structures on degree 27 extensions}

\begin{center}
\begin{tabular}{|c||c|c|c|c|c||c|}
\hline
\multicolumn{1}{|c||}{} & \multicolumn{6}{|c|}{\bf Hopf Galois structures} \\
\hline \hline
{\bf G/Type } &  \, $C_{27}$\,   & $C_9\times C_3$ &  \, $H_{27}$\,  &  \, $C_3^3$ \,  &  \, $G_{27}$\,   & Total\\
\hline
\hline
\multicolumn{1}{|c||}{$C_{27}$}         & 9 &0&0&0& 0 &9   \\ \hline
\multicolumn{1}{|c||}{$C_9\times C_3$} & 0 &39&12&6 & 78  &135   \\ \hline
\multicolumn{1}{|c||}{$H_{27}$} & 0&48&318&51&96 &513  \\ \hline
\multicolumn{1}{|c||}{$C_3^3$} & 0 &624&1326&339 & 1248 &3537   \\ \hline
\multicolumn{1}{|c||}{$G_{27}$}      &0&39&12  & 6  &78&135  \\ \hline
\multicolumn{1}{|c||}{$27T6$}       & 0 &0&78&27   & 0 &105\\ \hline
\multicolumn{1}{|c||}{$27T7$}                     & 0&0&0&1    & 0 &1\\ \hline
\multicolumn{1}{|c||}{$27T8$}                     & 1 &0&0&0    & 0 &1   \\ \hline
\multicolumn{1}{|c||}{$27T9$}                     & 0 &7&4&2   & 14&27  \\ \hline
\multicolumn{1}{|c||}{$27T10$}                     & 0 &1&0&0   & 0 &1  \\ \hline
\multicolumn{1}{|c||}{$27T11$}                     & 0 &4&22&5   & 8  &39  \\ \hline
\multicolumn{1}{|c||}{$27T12$}                     & 0&9&0&0     & 0&9 \\ \hline
\multicolumn{1}{|c||}{$27T13$}                     & 0 &16&94&35    & 32  &177  \\ \hline
\multicolumn{1}{|c||}{$27T14$}                     & 0 &7&4&2  & 14&27    \\ \hline
\multicolumn{1}{|c||}{$27T15$}                     & 0 &0&0&33    & 0&33   \\ \hline
\multicolumn{1}{|c||}{$27T16$}                     & 0&39&12&6   & 78 &135   \\ \hline
\multicolumn{1}{|c||}{$27T17$}                     & 0&9&0&0     & 18&27    \\ \hline
\multicolumn{1}{|c||}{$27T18$}                     & 0&12&120&33   & 24 &189 \\ \hline
\multicolumn{1}{|c||}{$27T19$}                     & 0&12&12&6   & 24 &54 \\ \hline
\multicolumn{1}{|c||}{$27T20$}                     & 0&9&0&0   & 18 &27 \\ \hline

\end{tabular}
\hspace{1cm}
\begin{tabular}{|c||c|c|c|c|c||c|}
\hline
\multicolumn{1}{|c||}{} & \multicolumn{6}{|c|}{\bf Hopf Galois structures} \\
\hline \hline
{\bf G/Type } &  \, $C_{27}$ \,  & $C_9\times C_3$ &  \, $H_{27}$\,  &  \, $C_3^3$ \,  &  \, $G_{27}$ \,  & Total\\
\hline
\hline
\multicolumn{1}{|c||}{$27T21$}         & 0 &6&6&3& 12 &27   \\ \hline
\multicolumn{1}{|c||}{$27T22$} & 9 &0&0&0 &0  &9  \\ \hline
\multicolumn{1}{|c||}{$27T23$} & 0&3&12&6&6 &27  \\ \hline
\multicolumn{1}{|c||}{$27T27$} & 0 &6&6&3 & 12 &27   \\ \hline
\multicolumn{1}{|c||}{$27T28$}      &0&27&0  & 0  &54&81  \\ \hline
\multicolumn{1}{|c||}{$27T29$}       & 0 &0&10&3  & 0 &13\\ \hline
\multicolumn{1}{|c||}{$27T30$}      & 0&1&0&0    & 0 &1\\ \hline
\multicolumn{1}{|c||}{$27T31$}                     & 0 &0&0&9   & 0 &9   \\ \hline
\multicolumn{1}{|c||}{$27T32$}                     & 0 &0&6&3   & 0&9 \\ \hline
\multicolumn{1}{|c||}{$27T33$}                     & 0 &0&6&3   & 0 &9  \\ \hline
\multicolumn{1}{|c||}{$27T34$}                     & 0 &0&0&1  & 0  &1  \\ \hline
\multicolumn{1}{|c||}{$27T35$}                     & 0&0&18&7     & 0&25 \\ \hline
\multicolumn{1}{|c||}{$27T36$}                     & 0 &0&0&9    & 0  &9  \\ \hline
\multicolumn{1}{|c||}{$27T37$}                     & 0 &0&0&6  & 0&6   \\ \hline
\multicolumn{1}{|c||}{$27T39$}                     & 0 &9&0&0    & 0&9   \\ \hline
\multicolumn{1}{|c||}{$27T46$}                     & 0&0&0&15   & 0 &15   \\ \hline
\multicolumn{1}{|c||}{$27T47$}                     & 0&3&0&0     & 0&3   \\ \hline
\multicolumn{1}{|c||}{$27T48$}                     & 0&0&0&3   & 0 &3 \\ \hline
\multicolumn{1}{|c||}{$27T49$}                     & 0&3&0&0  & 0 &3 \\ \hline
\multicolumn{1}{|c||}{$27T50$}                     & 0&0&6&3   & 0 &9 \\ \hline

\end{tabular}
\end{center}

\end{landscape}

\end{document}